\documentclass[12pt]{amsart}
\usepackage{amsfonts,amsmath,amsthm,amssymb}
\usepackage{latexsym}
\usepackage{paralist}
\newtheorem{theorem}{Theorem}[section]
\newtheorem{corollary}[theorem]{Corollary}

\newtheorem{lemma}[theorem]{Lemma}

\newtheorem{proposition}[theorem]{Proposition}

\theoremstyle{definition}
\theoremstyle{remark}

\theoremstyle{remark}

\newcommand{\beql}[1]{\begin{equation}\label{#1}}
\newcommand{\eeq}{\end{equation}}

\begin{document}
\title{A further look into combinatorial orthogonality}
\author{Simone Severini, Ferenc Sz\"oll\H{o}si}

\address{Simone Severini: Institute for Quantum Computing and
Department of Combinatorics and Optimization, University of
Waterloo, N2L 3G1 Waterloo, Canada} \email{simoseve@gmail.com}
\address{Ferenc Sz\"oll\H{o}si: Budapest University of Technology and Economics (BUTE), Institute of
Mathematics, Department of Analysis, Egry J. u. 1, H-1111 Budapest, Hungary}
\email{szoferi@math.bme.hu}

\begin{abstract}
Strongly quadrangular matrices have been introduced in the study of
the combinatorial properties of unitary matrices. It is known that
if a $(0,1)$-matrix supports a unitary then it is strongly
quadrangular. However, the converse is not necessarily true. In this
paper, we fully classify strongly quadrangular matrices up to degree
$5$. We prove that the smallest strongly quadrangular matrices which
do not support unitaries have exactly degree $5$. Further, we
isolate two submatrices not allowing a $(0,1)$-matrix to support
unitaries.
\end{abstract}

\maketitle

%\thanks{}



\section{Introduction}

Orthogonality is a common concept which generalizes perpendicularity
in Euclidean geometry. It appears in different mathematical
contexts, like linear algebra, functional analysis, combinatorics,
\emph{etc.} The necessary ingredient for introducing orthogonality
is a notion that allow to measure the \emph{angle} between two
objects. For example, in sufficiently rich vector spaces, this
consists of the usual inner product $\langle \cdot ,\cdot \rangle $.
Specifically, two vectors $u=(u_{1},..,u_{n})$ and
$v=(v_{1},...v_{n})$,
from a vector space over a generic field $\mathbb{F}$, are said to be \emph{%
orthogonal}, if $\langle u,v\rangle =0$. It is immediately clear that $u$
and $v$ are orthogonal only if there is a special relation between their
entries, and that this relation does not only involve the magnitude and the
signs, but also the position of the zeros, if there are any.

At a basic level, dealing with orthogonality from the combinatorial
point of view means, among other things, to study \emph{pattern} of
zeros in arrangements of vectors, some of which are orthogonal to
each other. A natural, somehow extremal scenario, is when the
vectors form a square matrix. Indeed, a matrix $M$ with entries on
$\mathbb{F}$ is said to be \emph{orthogonal} if $\langle
M_{i},M_{j}\rangle =0$, for every two different rows and columns,
$M_{i}$ and $M_{j}$. In this setting one can state the following
natural problem: characterize the zero pattern of orthogonal
matrices. This is a typical problem in \emph{combinatorial matrix
theory}, the field of matrix theory concerning intrinsic properties
of matrices viewed as arrays of numbers rather than algebraic
objects in themselves (see \cite{BR91}).

The term \textquotedblleft zero pattern\textquotedblright\ is not
followed by a neat mathematical definition. By zero pattern, we
intuitively mean the position of the zeros seen as forming a whole.
A more concrete definition may be introduced in the language of
graph theory. Let $D=(V,E)$ be a directed graph, without multiple
edges, but possibly with self-loops (see \cite{BG01} for this
standard graph theoretic terminology). Let $A(D)$ be the adjacency
matrix of $D$. We say that a matrix $M$ with entries on $\mathbb{F}
$ is \emph{supported} by $D$ (or, equivalently, by $A(D)$), if we
obtain $A(D)$ when replacing
with ones the nonzero entries of $M$. In other words, $D$ supports $M$, if $A(D)$ and $%
M$ have the same zero pattern. If this is the case, then $D$ is said
to be the \emph{digraph} \emph{of }$M$. Equivalently, $A(D)$ is said
to be the \emph{support} of $M$.

Studying the zero pattern of a family of matrices with certain
properties is equivalent to characterize the class of digraphs of
the matrices. When the field is $\mathbb{R}$ or $\mathbb{C}$, an
orthogonal matrix is also said to be \emph{real orthogonal} or
\emph{unitary}, respectively. These are practically ubiquitous
matrices, with roles spanning from coding theory to signal
processing, and from industrial screening experiments to the quantum
mechanics of closed systems,\emph{etc}. Historically, the problem of
characterizing zero patterns of orthogonal matrices was firstly
formulated by Fiedler \cite{fie, fie1}, and it is contextually
related to the more general problem of characterizing
ortho-stochastic matrices.

Even if not explicitly, the same problem can also be found in some
foundational issue of quantum theory (see \cite{lou}). Just
recently, this was motivated by some other questions concerning
unitary quantum evolution on graphs \cite{argo, tob}. Like many
other situations involving orthogonality, characterizing the zero
pattern of orthogonal matrices is not a simple problem. In some way,
a justification comes form the difficulty that we encounter when
trying to classify weighing, real and complex Hadamard matrices
\cite{MRSZ,Sz,karol}, and the related combinatorial designs
\cite{GS78} (see \cite{Ha}, for a more recent survey). One major
obstacle is in the \emph{global} features of orthogonality. Loosely
speaking, it is in fact evident that the essence of orthogonality
can not be isolated by looking at forbidden submatrices only, but
the property is subtler because it asks for relations between the
submatrices.

A first simple condition for orthogonality was proposed by Beasley, Brualdi
and Shader in 1991 \cite{BBS93}. As a tool, the authors introduced \emph{%
combinatorial orthogonality}. A $(0,1)$\emph{-matrix} is a matrix with
entries in the set $\{0,1\}$. A $(0,1)$-matrix $M$ is said to be \emph{%
combinatorially orthogonal}, or, equivalently, \emph{quadrangular},
if $\langle M_{i},M_{j}\rangle \neq 1$, for every two different rows
and columns, $M_{i}$ and $M_{j}$. It is immediate to observe that
the adjacency matrix of the digraph of an orthogonal matrix needs to
be combinatorially orthogonal. However, as it was already pointed
out in \cite{BBS93}, this condition is not sufficient to
characterize the zero pattern of orthogonal matrices. For the next
ten years, the few sporadic papers on this subject did focus on
quantitative results, mainly about the possible number of zeros
\cite{CJLP99,CS99,CS00a,CS00b}. In \cite{SE1} the problem was
reconsidered with the idea of pursuing a systematic study of the
qualitative side. The first step consisted of defining an easy
generalization of combinatorial orthogonality. This led to the
notion of \emph{strong quadrangularity}. Let $M$ be a (0,1)-matrix
of degree $n$, and let $S$ be a set of rows of $M$, forming an
$|S|\times n$ matrix. Suppose that for every $u\in S$ there exists a
row $v\in S$ such that $\langle u,v \rangle \neq 0$. Thus, if the
number of columns in $S$ containing at least two ones is at least
$|S|$ then $M$ is said to be \emph{row-strongly-quadrangular}. If
both $M$ and its transpose are row-strongly quadrangular then $M$ is
said to be \emph{strongly quadrangular} (for short, SQ). Even if
strong quadrangularity helps in exactly characterizing some classes
of digraphs of orthogonal matrices \cite{SE1}, the condition is not
necessary and sufficient. A counterexample involving a tournament
matrix of order $15$ was exhibited by Lundgren \emph{et al.}
\cite{B1}.

Let us denote by $\mathcal{U}_{n}$ the set of all (0,1)-matrices
whose digraph supports unitaries. Recall that an $n \times n$ matrix
$M$ is said to be \emph{indecomposable} if it has no $%
r\times (n-r)$ zero submatrix. The goal of the paper is to
investigate SQ matrices of small degree and find certain forbidden
substructures which prevent a (0,1)-matrix to support unitary
matrices.

One of the tools used through the paper is a construction due to
Di\c{t}\u{a} \cite{dita}, which is a generalization of the Kronecker
product. Although the original construction was defined for complex
Hadamard matrices, it can be easily extended to any unitary of
composite degree.

\begin{lemma}[Di\c{t}\u{a}'s construction]
Let $U_1,U_2,\hdots, U_{k}$ be unitaries of degree $n$, and let $%
[H]_{ij}=h_{ij}$ be a unitary of degree $k$. Then the following
matrix $Q$ of degree $nk$ is also unitary.

\begin{equation*}
Q=\left[%
\begin{array}{cccc}
h_{11}U_1 & h_{12}U_2 & \hdots & h_{1k}U_k \\
h_{21}U_1 & h_{22}U_2 & \hdots & h_{2k}U_k \\
\vdots & \vdots & \ddots & \vdots \\
h_{k1}U_1 & h_{k2}U_2 & \hdots & h_{kk}U_k \\
&  &  &
\end{array}%
\right],
\end{equation*}
\end{lemma}

\begin{corollary}
If $M_1, M_2,\hdots M_k \in \mathcal{U}_{n}$ then also the following
matrix $K \in \mathcal{U}_{kn}$:

\begin{equation*}
K= \left[%
\begin{array}{cccc}
M_1 & M_2 & \hdots & M_k \\
M_1 & M_2 & \hdots & M_k \\
\vdots & \vdots & \ddots & \vdots \\
M_1 & M_2 & \hdots & M_k \\
\end{array}%
\right].
\end{equation*}
\end{corollary}

\begin{proof}
Choose $H=\frac{1}{\sqrt{k}}F_k$, where $F_k$ is the matrix of the
Fourier transform over $\mathbb{Z}_{k}$, the abelian group of the
integers modulo $k$.
\end{proof}

It is clear why Di\c{t}\u{a}'s construction is useful for our
purposes. For example, the following matrix

\begin{equation*}
M=\left[%
\begin{array}{cccc}
1 & 0 & 1 & 1 \\
0 & 1 & 1 & 1 \\
1 & 0 & 1 & 1 \\
0 & 1 & 1 & 1 \\
\end{array}%
\right],
\end{equation*}

is in $\mathcal{U}_{4}$, since we can choose $U_1=I_2,
U_2=\frac{1}{\sqrt{2}}F_2$ and then apply the construction.

\section{SQ matrices of small degree}\label{list}
The purpose of this section is twofold. On the one hand, we would
like to give a detailed list of SQ matrices of small degree. This is
done in the perspective of further investigation. On the other hand,
we directly enumerate indecomposable, SQ matrices up to degree $5$.
The method that we adopt in this enumeration is a three-step
procedure. First, we simply construct \emph{all} $(0,1)$-matrices of
degree $n\leq5$. Second, we exclude from this list matrices which
are not SQ or contains a line (row or column) of zeros. Finally, we
determine representative from equivalence classes of the remaining
matrices. We also compute the order of their automorphism group.
Recall that two matrices $M_{1}$ and $M_{2}$ are said to be
\emph{equivalent} if there are permutation matrices $P$ and $Q$
such that $PM_{1}Q=M_{2}$. As usual, the automorphism group of a
$(0,1)$-matrix $M$ is the set of ordered pairs $(P,Q)$ of
permutation matrices such that $PMQ=M$. Some heuristics helps to
simplify our task.

\begin{lemma}
Two $(0,1)$-matrices having a different number of zeros are not
equivalent.
\end{lemma}

Of course, the converse statement is not necessarily true (see,
\emph{e.g.}, the two matrices of degree $4$ with exactly four zeros
below). Another natural heuristic consists of the number of ones in
each row. For a given (0,1)-matrix $[M]_{i,j}=m_{ij}$, we define the
multiset $\Lambda =\left\{ \sum_{j}m_{ij}:i=1,...,n\right\} $. Thus
the following observation is easy verify:

\begin{lemma}
Two $(0,1)$-matrices with different $\Lambda$'s are not equivalent.
\end{lemma}

Again, the converse statement does not hold in general. Recall that
a (0,1)-matrix is said to be \emph{regular} if the elements of
$\Lambda$ are all equal. The following lemma is specifically useful
for distinguishing regular matrices:

\begin{lemma}
Two $(0,1)$-matrices with nonisomorphic automorphism group are not
equivalent.
\end{lemma}

It follows that two $(0,1)$-matrices with automorphism groups of
different order are not equivalent. Unfortunately, there are
examples of nonequivalent matrices whose automorphism groups are
isomorphic. In particular, it
might happen that a matrix and its transpose are not equivalent. By
combining together the above facts, with the help of a computer, we
can fully classify matrices in $\mathcal{U}_{n}$ for $n \leq 5$. By
a careful analysis of the results, in Section \ref{beyond} we are
able to describe certain cases in which a matrix $M \notin
\mathcal{U}_{n}$ even if it is SQ. The smallest such an example is
of degree $5$. Additionally, if $M \in \mathcal{U}_{n}$ for $n< 5$
then $M$ is SQ and viceversa. We hereby present, up to equivalence,
the list of all indecomposable, SQ matrices of degree $n \leq 5$. We
need to fix some notational convention: if a matrix is equivalent to
a symmetric one we index it by $S$; if a matrix is not equivalent to
its transpose we index it by $T$. Regular matrices will be indexed
by $R$. Finally, the order of the automorphism group of a matrix is
written as a subscript. This information describes the number of
equivalent matrices in a given class. In particular,
$\left(n!\right)^2=\left|\mathrm{Aut}M\right|\cdot\#\left\{\text{Equivalent
matrices to $M$}\right\}.$

\subsection{n=1}

\beql{}
\left\{ \left[
\begin{array}{c}
1
\end{array}
\right]^{RS}_1\right\}
\eeq

This matrix, and more generally, every all-one matrix $J_n$, clearly
supports unitaries, since there is an $n \times n$ complex Hadamard
matrix for any $n$ \cite{Ha,karol}.

\subsection{n=2}

\beql{}
\left\{ \left[
\begin{array}{cc}
1 & 1 \\
1 & 1 \\
\end{array}
\right]^{RS}_4\right\}
\eeq

\subsection{n=3}

\beql{}
\left\{ \left[
\begin{array}{ccc}
0 & 1 & 1 \\
1 & 1 & 1 \\
1 & 1 & 1 \\
\end{array}
\right]^{S}_4, \left[
\begin{array}{ccc}
1 & 1 & 1 \\
1 & 1 & 1 \\
1 & 1 & 1 \\
\end{array}
\right]^{RS}_{36}\right\}
\eeq

\subsection{n=4}

\beql{}
\left\{ \left[
\begin{array}{cccc}
0 & 1 & 1 & 1 \\
1 & 1 & 1 & 1 \\
1 & 1 & 1 & 1 \\
1 & 1 & 1 & 1 \\
\end{array}
\right]^{S}_{36}, \left[
\begin{array}{cccc}
0 & 1 & 1 & 1 \\
1 & 0 & 1 & 1 \\
1 & 1 & 1 & 1 \\
1 & 1 & 1 & 1 \\
\end{array}
\right]^{S}_{8}, \left[
\begin{array}{cccc}
0 & 1 & 1 & 1 \\
1 & 0 & 1 & 1 \\
1 & 1 & 0 & 1 \\
1 & 1 & 1 & 1 \\
\end{array}
\right]^{S}_{6}, \left[
\begin{array}{cccc}
0 & 1 & 1 & 1 \\
1 & 0 & 1 & 1 \\
1 & 1 & 0 & 1 \\
1 & 1 & 1 & 0 \\
\end{array}
\right]^{RS}_{24},\right.
\eeq

\begin{equation*}
\left. \left[
\begin{array}{cccc}
0 & 0 & 1 & 1 \\
0 & 1 & 1 & 1 \\
1 & 1 & 1 & 1 \\
1 & 1 & 1 & 1 \\
\end{array}
\right]^{S}_{4}, \left[
\begin{array}{cccc}
0 & 0 & 1 & 1 \\
1 & 1 & 1 & 1 \\
1 & 1 & 1 & 1 \\
1 & 1 & 1 & 1 \\
\end{array}
\right]^{T}_{24}, \left[
\begin{array}{cccc}
1 & 0 & 1 & 1 \\
0 & 1 & 1 & 1 \\
1 & 0 & 1 & 1 \\
0 & 1 & 1 & 1 \\
\end{array}
\right]^{T}_{16}, \left[
\begin{array}{cccc}
1 & 1 & 1 & 1 \\
1 & 1 & 1 & 1 \\
1 & 1 & 1 & 1 \\
1 & 1 & 1 & 1 \\
\end{array}
\right]^{RS}_{576} \right\}
\end{equation*}

With the data above and going through all the few decomposable
matrices, we can give the following statement:

\begin{proposition} A $(0,1)$-matrix of degree $n\leq 4$ supports a unitary if and only
if it is SQ.
\end{proposition}

We will see later that this is not in general the case.

\subsection{n=5}

The following list contains $63$ items. Here, we double count the
matrices with index $T$. We can observe that not all of these
support unitaries, as we will see in Section \ref{beyond}.

\scriptsize{
%1-4
\beql{} \left\{\left[\begin{array}{ccccc}
 0 & 1 & 1 & 1 & 1 \\
 1 & 1 & 1 & 1 & 1 \\
 1 & 1 & 1 & 1 & 1 \\
 1 & 1 & 1 & 1 & 1 \\
 1 & 1 & 1 & 1 & 1
\end{array}
\right]^{S}_{576},\left[
\begin{array}{ccccc}
 0 & 1 & 1 & 1 & 1 \\
 1 & 0 & 1 & 1 & 1 \\
 1 & 1 & 1 & 1 & 1 \\
 1 & 1 & 1 & 1 & 1 \\
 1 & 1 & 1 & 1 & 1
\end{array}
\right]^{S}_{72},\left[
\begin{array}{ccccc}
 0 & 1 & 1 & 1 & 1 \\
 1 & 0 & 1 & 1 & 1 \\
 1 & 1 & 0 & 1 & 1 \\
 1 & 1 & 1 & 1 & 1 \\
 1 & 1 & 1 & 1 & 1
\end{array}
\right]^{S}_{24},\left[
\begin{array}{ccccc}
 0 & 1 & 1 & 1 & 1 \\
 1 & 0 & 1 & 1 & 1 \\
 1 & 1 & 0 & 1 & 1 \\
 1 & 1 & 1 & 0 & 1 \\
 1 & 1 & 1 & 1 & 1
\end{array}
\right]^{S}_{24},\right.
\eeq
%5-8
\[\left[
\begin{array}{ccccc}
 0 & 1 & 1 & 1 & 1 \\
 1 & 0 & 1 & 1 & 1 \\
 1 & 1 & 0 & 1 & 1 \\
 1 & 1 & 1 & 0 & 1 \\
 1 & 1 & 1 & 1 & 0
\end{array}
\right]^{RS}_{120},
\left[
\begin{array}{ccccc}
 0 & 0 & 0 & 1 & 1 \\
 0 & 0 & 1 & 1 & 1 \\
 0 & 1 & 1 & 1 & 1 \\
 1 & 1 & 1 & 1 & 1 \\
 1 & 1 & 1 & 1 & 1
\end{array}
\right]^{S}_{4},
\left[
\begin{array}{ccccc}
 0 & 0 & 0 & 1 & 1 \\
 0 & 1 & 1 & 1 & 1 \\
 0 & 1 & 1 & 1 & 1 \\
 1 & 1 & 1 & 1 & 1 \\
 1 & 1 & 1 & 1 & 1
\end{array}
\right]^{S}_{16},\left[
\begin{array}{ccccc}
 0 & 0 & 0 & 1 & 1 \\
 0 & 1 & 1 & 1 & 1 \\
 0 & 1 & 1 & 1 & 1 \\
 1 & 1 & 1 & 0 & 0 \\
 1 & 1 & 1 & 0 & 0
\end{array}
\right]^{S}_{16},\]
%9-12
\[
\left[
\begin{array}{ccccc}
 0 & 0 & 1 & 1 & 1 \\
 0 & 1 & 1 & 1 & 1 \\
 1 & 1 & 1 & 1 & 1 \\
 1 & 1 & 1 & 1 & 1 \\
 1 & 1 & 1 & 1 & 1
\end{array}
\right]^{S}_{36},\left[
\begin{array}{ccccc}
 0 & 0 & 1 & 1 & 1 \\
 0 & 1 & 1 & 1 & 1 \\
 1 & 1 & 0 & 1 & 1 \\
 1 & 1 & 1 & 1 & 1 \\
 1 & 1 & 1 & 1 & 1
\end{array}
\right]^{S}_{4},
\left[
\begin{array}{ccccc}
 0 & 0 & 1 & 1 & 1 \\
 0 & 1 & 1 & 1 & 1 \\
 1 & 1 & 0 & 1 & 1 \\
 1 & 1 & 1 & 0 & 1 \\
 1 & 1 & 1 & 1 & 1
\end{array}
\right]^{S}_{2},
\left[
\begin{array}{ccccc}
 0 & 0 & 1 & 1 & 1 \\
 0 & 1 & 1 & 1 & 1 \\
 1 & 1 & 0 & 1 & 1 \\
 1 & 1 & 1 & 0 & 1 \\
 1 & 1 & 1 & 1 & 0
\end{array}
\right]^{S}_{6},\]
%13-16
\[\left[
\begin{array}{ccccc}
 0 & 0 & 1 & 1 & 1 \\
 0 & 0 & 1 & 1 & 1 \\
 1 & 1 & 1 & 1 & 1 \\
 1 & 1 & 1 & 1 & 1 \\
 1 & 1 & 1 & 1 & 1
\end{array}
\right]^{S}_{144},
\left[
\begin{array}{ccccc}
 0 & 0 & 1 & 1 & 1 \\
 0 & 0 & 1 & 1 & 1 \\
 1 & 1 & 0 & 1 & 1 \\
 1 & 1 & 1 & 1 & 1 \\
 1 & 1 & 1 & 1 & 1
\end{array}
\right]^{S}_{16},
\left[
\begin{array}{ccccc}
 0 & 0 & 1 & 1 & 1 \\
 0 & 0 & 1 & 1 & 1 \\
 1 & 1 & 0 & 1 & 1 \\
 1 & 1 & 1 & 0 & 1 \\
 1 & 1 & 1 & 1 & 1
\end{array}
\right]^{NS}_{8},\left[
\begin{array}{ccccc}
 0 & 0 & 1 & 1 & 1 \\
 0 & 0 & 1 & 1 & 1 \\
 1 & 1 & 0 & 1 & 1 \\
 1 & 1 & 1 & 0 & 1 \\
 1 & 1 & 1 & 1 & 0
\end{array}
\right]^{NS}_{24},
\]
%17-20
\[
\left[
\begin{array}{ccccc}
 0 & 0 & 1 & 1 & 1 \\
 0 & 1 & 0 & 1 & 1 \\
 1 & 0 & 1 & 1 & 1 \\
 1 & 1 & 1 & 1 & 1 \\
 1 & 1 & 1 & 1 & 1
\end{array}
\right]^{S}_4,
\left[
\begin{array}{ccccc}
 0 & 0 & 1 & 1 & 1 \\
 0 & 1 & 0 & 1 & 1 \\
 1 & 0 & 1 & 1 & 1 \\
 1 & 1 & 1 & 0 & 1 \\
 1 & 1 & 1 & 1 & 1
\end{array}
\right]^{S}_1,
\left[
\begin{array}{ccccc}
 0 & 0 & 1 & 1 & 1 \\
 0 & 1 & 0 & 1 & 1 \\
 1 & 0 & 1 & 1 & 1 \\
 1 & 1 & 1 & 0 & 1 \\
 1 & 1 & 1 & 1 & 0
\end{array}
\right]^{S}_2,
\left[
\begin{array}{ccccc}
 1 & 0 & 0 & 1 & 1 \\
 0 & 1 & 1 & 1 & 1 \\
 0 & 1 & 1 & 1 & 1 \\
 1 & 1 & 1 & 1 & 1 \\
 1 & 1 & 1 & 1 & 1
\end{array}
\right]^{S}_{16},\]
%21-24
\[
\left[
\begin{array}{ccccc}
 1 & 0 & 0 & 1 & 1 \\
 0 & 1 & 1 & 1 & 1 \\
 0 & 1 & 1 & 1 & 1 \\
 1 & 1 & 1 & 0 & 1 \\
 1 & 1 & 1 & 1 & 1
\end{array}
\right]^{S}_{4},
\left[
\begin{array}{ccccc}
 1 & 0 & 0 & 1 & 1 \\
 0 & 1 & 1 & 1 & 1 \\
 0 & 1 & 1 & 1 & 1 \\
 1 & 1 & 1 & 0 & 1 \\
 1 & 1 & 1 & 1 & 0
\end{array}
\right]^{S}_{8},
\left[
\begin{array}{ccccc}
 1 & 1 & 1 & 1 & 1 \\
 1 & 1 & 1 & 1 & 1 \\
 1 & 1 & 1 & 1 & 1 \\
 1 & 1 & 1 & 1 & 1 \\
 1 & 1 & 1 & 1 & 1
\end{array}
\right]^{RS}_{14400},
\left[
\begin{array}{ccccc}
 0 & 0 & 1 & 1 & 1 \\
 1 & 1 & 1 & 1 & 1 \\
 1 & 1 & 1 & 1 & 1 \\
 1 & 1 & 1 & 1 & 1 \\
 1 & 1 & 1 & 1 & 1
\end{array}
\right]^{T}_{288},\]
%25-28
\[
\left[
\begin{array}{ccccc}
 0 & 0 & 0 & 1 & 1 \\
 1 & 1 & 1 & 1 & 1 \\
 1 & 1 & 1 & 1 & 1 \\
 1 & 1 & 1 & 1 & 1 \\
 1 & 1 & 1 & 1 & 1
\end{array}
\right]^{T}_{288},\left[
\begin{array}{ccccc}
 0 & 0 & 0 & 1 & 1 \\
 0 & 1 & 1 & 1 & 1 \\
 1 & 1 & 1 & 1 & 1 \\
 1 & 1 & 1 & 1 & 1 \\
 1 & 1 & 1 & 1 & 1
\end{array}
\right]^{T}_{24},\left[
\begin{array}{ccccc}
 0 & 0 & 0 & 1 & 1 \\
 0 & 1 & 1 & 1 & 1 \\
 1 & 0 & 1 & 1 & 1 \\
 1 & 1 & 1 & 1 & 1 \\
 1 & 1 & 1 & 1 & 1
\end{array}
\right]^{T}_{8},\left[
\begin{array}{ccccc}
 0 & 0 & 0 & 1 & 1 \\
 0 & 1 & 1 & 1 & 1 \\
 1 & 0 & 1 & 1 & 1 \\
 1 & 1 & 0 & 1 & 1 \\
 1 & 1 & 1 & 1 & 1
\end{array}
\right]^{T}_{12},\]
%28-32
\[
\left[
\begin{array}{ccccc}
 0 & 0 & 0 & 1 & 1 \\
 1 & 1 & 1 & 1 & 1 \\
 1 & 1 & 1 & 1 & 1 \\
 1 & 1 & 1 & 1 & 1 \\
 1 & 1 & 1 & 0 & 0
\end{array}
\right]^{T}_{72},\left[
\begin{array}{ccccc}
 0 & 0 & 0 & 1 & 1 \\
 0 & 1 & 1 & 1 & 1 \\
 1 & 1 & 1 & 1 & 1 \\
 1 & 1 & 1 & 1 & 1 \\
 1 & 1 & 1 & 0 & 0
\end{array}
\right]^{T}_{8},\left[
\begin{array}{ccccc}
 0 & 0 & 0 & 1 & 1 \\
 0 & 1 & 1 & 1 & 1 \\
 1 & 0 & 1 & 1 & 1 \\
 1 & 1 & 1 & 1 & 1 \\
 1 & 1 & 1 & 0 & 0
\end{array}
\right]^{T}_{4},
\left[
\begin{array}{ccccc}
 0 & 0 & 0 & 1 & 1 \\
 0 & 1 & 1 & 1 & 1 \\
 1 & 0 & 1 & 1 & 1 \\
 1 & 1 & 0 & 1 & 1 \\
 1 & 1 & 1 & 0 & 0
\end{array}
\right]^{T}_{12},\]
%22-36
\[
\left[
\begin{array}{ccccc}
 0 & 0 & 1 & 1 & 1 \\
 0 & 1 & 0 & 1 & 1 \\
 1 & 1 & 1 & 1 & 1 \\
 1 & 1 & 1 & 1 & 1 \\
 1 & 1 & 1 & 1 & 1
\end{array}
\right]^{T}_{24},\left[
\begin{array}{ccccc}
 0 & 0 & 1 & 1 & 1 \\
 0 & 1 & 0 & 1 & 1 \\
 1 & 1 & 1 & 0 & 1 \\
 1 & 1 & 1 & 1 & 1 \\
 1 & 1 & 1 & 1 & 1
\end{array}
\right]^{T}_{4},
\left[
\begin{array}{ccccc}
 0 & 0 & 1 & 1 & 1 \\
 0 & 1 & 0 & 1 & 1 \\
 1 & 1 & 1 & 0 & 1 \\
 1 & 1 & 1 & 1 & 0 \\
 1 & 1 & 1 & 1 & 1
\end{array}
\right]^{T}_{4},\left[
\begin{array}{ccccc}
 0 & 0 & 0 & 1 & 1 \\
 0 & 1 & 1 & 0 & 0 \\
 1 & 1 & 1 & 1 & 1 \\
 1 & 1 & 1 & 1 & 1 \\
 1 & 1 & 1 & 1 & 1
\end{array}
\right]^{T}_{48},\]
%36-40
\[
\left[
\begin{array}{ccccc}
 0 & 0 & 1 & 1 & 1 \\
 1 & 1 & 0 & 1 & 1 \\
 1 & 1 & 1 & 1 & 1 \\
 1 & 1 & 1 & 1 & 1 \\
 1 & 1 & 1 & 1 & 1
\end{array}
\right]^{T}_{24},\left[
\begin{array}{ccccc}
 0 & 0 & 1 & 1 & 1 \\
 1 & 1 & 0 & 1 & 1 \\
 1 & 1 & 1 & 0 & 1 \\
 1 & 1 & 1 & 1 & 1 \\
 1 & 1 & 1 & 1 & 1
\end{array}
\right]^{T}_{8},
\left[
\begin{array}{ccccc}
 0 & 0 & 1 & 1 & 1 \\
 1 & 1 & 0 & 1 & 1 \\
 1 & 1 & 1 & 0 & 1 \\
 1 & 1 & 1 & 1 & 0 \\
 1 & 1 & 1 & 1 & 1
\end{array}
\right]^{T}_{12},
\left[
\begin{array}{ccccc}
 0 & 0 & 0 & 1 & 1 \\
 1 & 1 & 1 & 0 & 0 \\
 1 & 1 & 1 & 0 & 0 \\
 1 & 1 & 1 & 1 & 1 \\
 1 & 1 & 1 & 1 & 1
\end{array}
\right]^{T}_{48},\]
%41-43
\[
\left.
\left[
\begin{array}{ccccc}
 0 & 0 & 0 & 1 & 1 \\
 0 & 0 & 1 & 1 & 1 \\
 1 & 1 & 1 & 1 & 1 \\
 1 & 1 & 1 & 1 & 1 \\
 1 & 1 & 1 & 1 & 1
\end{array}
\right]^{T}_{24},
\left[
\begin{array}{ccccc}
 0 & 0 & 0 & 1 & 1 \\
 0 & 0 & 1 & 1 & 1 \\
 1 & 1 & 0 & 0 & 0 \\
 1 & 1 & 1 & 1 & 1 \\
 1 & 1 & 1 & 1 & 1
\end{array}
\right]^{T}_{8},\left[
\begin{array}{ccccc}
 0 & 0 & 0 & 1 & 1 \\
 0 & 1 & 1 & 0 & 0 \\
 0 & 1 & 1 & 1 & 1 \\
 1 & 1 & 1 & 1 & 1 \\
 1 & 1 & 1 & 1 & 1
\end{array}
\right]^{T}_{16}\right\}\]

%tiny
}
\normalsize

While constructing unitaries matching a given pattern up to degree
$4$ is a simple task, considering $n=5$ brings up several
difficulties. First of all, as one can see, there are many
equivalent classes and presenting unitaries for each and every class
is out of reach. Secondly, it turns out that there are at least two
such matrices which do not support a unitary. We index these
matrices by $N$. This statement will be formally proved in Section
\ref{beyond}.

Since the number of SQ matrices grows very fast, lacking of
computational power, we did stop our counting at $n=5$. However, we
propose two further special cases which are arguably easier to
handle.

\subsubsection{Symmetric SQ matrices}

It is evident that the main difficulty in classifying SQ matrices is
\emph{not} the actual construction of the matrices, but determining
equivalence classes. This is a time-consuming procedure even for
small degrees. The following lemma shows that classifying only
\emph{symmetric} SQ matrices is a definitely easier problem.

\begin{lemma}\label{sym}
If a $(0,1)$-matrix $M$ is equivalent to a symmetric one, then there
is a permutation matrix $R$, such that $RMR=M^T$.
\end{lemma}

\begin{proof}
Suppose that $M$ is equivalent to a symmetric matrix, denoted by
$S$. Then there are permutation matrices $P$ and $Q$, such that
$PMQ=S=S^T$, so $PMQ=Q^TM^TP^T$, and hence
$\left(QP\right)M\left(QP\right)=M^T$. This implies that $R=QP$ is a
permutation matrix, as required.
\end{proof}

Determining whether a $(0,1)$-matrix $M$ is equivalent to a symmetric
one therefore simply boils down to a two phase procedure: first, we
check if there are permutations matrices for which $RMR=M^T$;
second, we check if $Q^TRMQ$ is symmetric for a certain $Q$. If
there exists such a pair of permutation matrices $R$ and $Q$, then
$M$ is equivalent to a symmetric matrix. This procedure is clearly
faster than simultaneously looking for $P$ and $Q$ such that $PMQ$
is symmetric.

\subsubsection{Regular SQ matrices}

Here we focus on \emph{regular} SQ matrices. We have classified
these matrices up to degree $6$. The results up to degree $5$ can be
found in the lists above. The list for degree $6$ is included below.
Let $\sigma$ be the number of nonzero entries in each row of a
regular matrix. There are regular SQ matrices of order $6$ with
$\sigma=6,5,3,2,1$, since $J_6, J_6-I_6, I_3\oplus I_3, I_2\oplus
I_2\oplus I_2, I_6$ are such examples, where $I_n$ denotes the
$n\times n$ identity matrix. It can be checked that in fact these
are the only ones. However, the case $\sigma=4$ turns out to be
interesting, since one out of the four regular matrices does not
support unitaries. This fact will be investigated later in Theorem
\ref{cond} of Section \ref{beyond}.

\scriptsize{
\beql{}
\left\{
\left[
\begin{array}{cccccc}
 1 & 0 & 0 & 1 & 1 & 1 \\
 0 & 1 & 0 & 1 & 1 & 1 \\
 0 & 0 & 1 & 1 & 1 & 1 \\
 1 & 1 & 1 & 1 & 0 & 0 \\
 1 & 1 & 1 & 0 & 1 & 0 \\
 1 & 1 & 1 & 0 & 0 & 1
\end{array}
\right]^{RS}_{72},
\left[
\begin{array}{cccccc}
 0 & 0 & 1 & 1 & 1 & 1 \\
 0 & 0 & 1 & 1 & 1 & 1 \\
 1 & 1 & 0 & 0 & 1 & 1 \\
 1 & 1 & 0 & 0 & 1 & 1 \\
 1 & 1 & 1 & 1 & 0 & 0 \\
 1 & 1 & 1 & 1 & 0 & 0
\end{array}
\right]^{RS}_{384},
\left[
\begin{array}{cccccc}
 1 & 1 & 1 & 1 & 0 & 0 \\
 0 & 1 & 1 & 1 & 1 & 0 \\
 0 & 0 & 1 & 1 & 1 & 1 \\
 1 & 0 & 0 & 1 & 1 & 1 \\
 1 & 1 & 0 & 0 & 1 & 1 \\
 1 & 1 & 1 & 0 & 0 & 1
\end{array}
\right]^{RS}_{12},
\left[
\begin{array}{cccccl}
 0 & 0 & 1 & 1 & 1 & 1 \\
 0 & 0 & 1 & 1 & 1 & 1 \\
 1 & 1 & 1 & 1 & 0 & 0 \\
 1 & 1 & 0 & 1 & 1 & 0 \\
 1 & 1 & 0 & 0 & 1 & 1 \\
 1 & 1 & 1 & 0 & 0 & 1
\end{array}
\right]^{NRS}_{32}\right\}
\eeq
}

\normalsize

We conclude by summarizing our observations:

\begin{itemize}
\item The number of inequivalent indecomposable SQ matrices of degree $%
n=1,2,...,5$ is $1,1,2,10,63$, respectively. All known terms of this
sequence match the number of triples of standard tableaux with the
same shape of height less than or equal to three. This sequence is
A129130 in \cite{sl}.

\item The number of inequivalent SQ matrices of orders $n=1,2,...,5$ is $%
1,2,4,15,80$, respectively.

\item The number of inequivalent indecomposable symmetric SQ matrices of
degree $n=1,2,...,5$ is $1,1,2,6,23$, respectively.

\item The number of inequivalent symmetric SQ matrices of
degree $n=1,2,...,5$ is\\ $1,2,4,11,44$, respectively.

\item The number of inequivalent indecomposable regular SQ matrices of
orders $n=1,2,...,6$ is $1,1,1,2,2,4$.

\item The number of inequivalent regular SQ matrices of degree $n=1,2,...,6$
is $1,2,2,4,3,9$, respectively.

\end{itemize}

\section{Beyond strong quadrangularity}\label{beyond}

In \cite{B1}, the authors exhibited the adjacency matrix of a
tournament on $15$ vertices, which, despite being SQ, it is not in
$\mathcal{U}_{15}$. The first result of this section is a refined
version of that. Specifically, we have the following:

\begin{theorem}{\label{cond}} Let $M$ (or its transpose) be a $(0,1)$-matrix equivalent
to a matrix in the following form, for $k\geq 1$:

\beql{}
\left[
\begin{array}{ccc}
Q & J_{3\times k} & X \\
Y & Z & \ast \\
\end{array}
\right],\ \mathnormal{ where }\ Q=\left[
\begin{array}{cc}
1 & 0 \\
0 & 1 \\
1 & 1 \\
\end{array}
\right].
\eeq

Further, suppose that
\begin{enumerate}
\item the rows of $X$ are mutually orthogonal,
\item every column of $Y$ is orthogonal to every column of $Z$.
\end{enumerate}
Then $M$ does not support unitaries.
\end{theorem}

\begin{proof}
The idea of the proof is exactly the same as in \cite{B1}. Suppose
on the contrary that there exists a unitary $U$ whose support is
$M$. Let $R_i$ and $C_i$ denote the $i$-th row and column of $U$
respectively, for each $i=1,\hdots n$ and let $[U]_{ij}=u_{i,j}$.
Now observe, that $\left\langle
C_1,C_j\right\rangle=u_{1,1}\overline{u}_{1,j}+u_{3,1}\overline{u}_{3,j}=0$,
where $j=3,4,\hdots,k+2$. This implies
$-\overline{u}_{1,1}/\overline{u}_{3,1}=u_{3,j}/u_{1,j}$, where
$j=3,4,\hdots,k+2$. So the vectors
$\left[u_{1,3},\hdots,u_{1,k+2}\right]$ and
$\left[u_{3,3},\hdots,u_{3,k+2}\right]$ are scalar multiples of each
other. Similarly: $\left\langle
C_2,C_j\right\rangle=u_{2,2}\overline{u}_{2,j}+u_{3,2}\overline{u}_{3,j}=0$,
where $j=3,4,\hdots,k+2$. So, this implies
$-\overline{u}_{2,2}/\overline{u}_{3,2}=u_{3,j}/u_{2,j}$, where
$j=3,4,\hdots,k+2$. So the vectors
$\left[u_{2,3},\hdots,u_{2,k+2}\right]$ and
$\left[u_{3,3},\hdots,u_{3,k+2}\right]$ are scalar multiples of each
other. It follows that $\left\langle R_1,R_2
\right\rangle=\left\langle
\left[u_{1,3},\hdots,u_{1,k+2}\right],\left[u_{2,3},\hdots,u_{2,k+2}\right]\right\rangle\neq0$,
a contradiction.
\end{proof}

The next statement summarizes the main features of the matrices satisfying the conditions of Theorem \ref{cond}.

\begin{proposition}\label{cond2}
Under the conditions of Theorem \ref{cond}, a SQ matrix of degree
$n$ satisfies the following properties:

\begin{enumerate}
\item $k\geq 2$;
\item The first row of $Y$ is $\left[1,1\right]$;
\item The first row of $Z$ is $\left[0,\hdots, 0\right]$;
\item $X$ has at least two columns;
\item $n\geq 6$.
\end{enumerate}
\end{proposition}

\begin{proof}
Suppose that we have a matrix equivalent to $M$. Since its first two
rows of share a common $1$, and $X$ cannot have two rows who share a
common $1$, $k\geq 2$ follows. Similarly, the first and second
column of $M$ share a common $1$, hence by quadrangularity, these
share another $1$, and up to equivalence, we can suppose that it is
in the $4$-th row of $M$. Thus, the first row of $Y$ can be chosen
to be $\left[1,1\right]$. By the second condition of Theorem
\ref{cond}, the first row of $Z$ should be $\left[0,\hdots,
0\right]$. Again, since the $1$-st and $4$-th rows share a
common $1$, by quadrangularity, they should share another $1$.
However, we have already seen that the first row of $Z$ is all $0$.
Thus, these rows must share this specific $1$ in $X$. The same
argument applies for the $2$-nd and $4$-th row of $M$, and since two
rows of $X$ cannot share a common $1$, it must have at leats two
columns. It follows that $k\leq n-4$ and therefore $n\geq 6$.
\end{proof}

Next, we estimate the possible number of $1$s in matrices satisfying
the conditions of Theorem \ref{cond}.

\begin{lemma}\label{count1}
Suppose that a SQ matrix $M$ of degree $n\geq 6$ satisfies the
conditions of Theorem \ref{cond}. Then, its possible number of ones
is at most $n^2-3n+6$, and hence, it has a least $3n-6$ zeros.
\end{lemma}

\begin{proof}
We simply count the number of ones in all blocks of $M$ separately.
First, the number of ones in $Q$ is $4$, and clearly, the number of
$1$s in $J$ are $3k$. By the first condition of Theorem \ref{cond},
the number of ones in $X$ is at most $n-k-2$. Now by Lemma
\ref{cond2} the first row of $Z$ is $\left[0,\hdots,0\right]$ (up to
equivalence), hence the number of ones in $Y$ and $Z$ together is at
most $k(n-4)+2$, and finally the number of ones in the lower right
submatrix is at most $(n-k-2)(n-3)$. Thus the possible number of
ones is \beql{} 4+3k+n-k-2+k(n-4)+2+(n-k-2)(n-3)=n^2-4n+10+k\leq
n^2-3n+6. \eeq
\end{proof}

The $6\times 6$ matrix $A$ below is SQ. However, by Theorem
\ref{cond}, $A \notin \mathcal{U}_{6}$.

\beql{}
A=\left[
\begin{array}{cc|cc|cc}
1 & 0 & 1 & 1 & 1 & 0 \\
0 & 1 & 1 & 1 & 0 & 1 \\
1 & 1 & 1 & 1 & 0 & 0 \\ \hline
1 & 1 & 0 & 0 & 1 & 1 \\
1 & 1 & 0 & 0 & 1 & 1 \\
0 & 0 & 1 & 1 & 1 & 1 \\
\end{array}
\right] \eeq

Note that $A$ is regular, therefore it is equivalent to the
exceptional regular matrix of degree $6$ appearing in Section
\ref{list}.

The example above shows that there are indeed SQ matrices of degree
$6$, which cannot support unitaries. It is of particular interest to
find out if there are such exceptional matrices already for degree
$5$. Lemma \ref{cond2} explains that we cannot rely on Theorem
\ref{cond}, since this result does not say anything about matrices
of order $5$. By analyzing the list of Section \ref{list}, one can
observe that such exceptional matrices do exist for degree $5$. The
reason for this phenomenon is summarized in the following:

\begin{theorem}\label{newcond}
Let $M$ (or its transpose) be a $(0,1)$-matrix equivalent to a matrix in the
following form:

\beql{}
\left[
\begin{array}{ccc}
Q & J_{3\times 2} & \ast \\
X & Y & \ast \\
\end{array}
\right],\ \mathnormal{where}\ Q=\left[
\begin{array}{cc}
1 & 0 \\
0 & 1 \\
1 & 1 \\
\end{array}
\right]
\eeq

Further, suppose that
\begin{enumerate}
\item the columns of $Y$ are mutually orthogonal,
\item every column of $X$ is orthogonal to every column of $Y$.
\end{enumerate}

Then $M$ does not support unitaries.
\end{theorem}

\begin{proof}
Suppose on the contrary that we have a unitary $U$, whose support is
$M$. Let us use the same notations as in the proof of Theorem
\ref{cond}. By orthogonality $\left\langle
C_1,C_3\right\rangle=u_{11}\overline{u}_{13}+u_{31}\overline{u}_{33}=0$,
$\left\langle
C_1,C_4\right\rangle=u_{11}\overline{u}_{14}+u_{31}\overline{u}_{34}=0$,
$\left\langle
C_2,C_3\right\rangle=u_{22}\overline{u}_{23}+u_{32}\overline{u}_{33}=0$,
$\left\langle
C_2,C_4\right\rangle=u_{22}\overline{u}_{24}+u_{32}\overline{u}_{34}=0$,
hence $u_{31}=-u_{11}\overline{u}_{13}/\overline{u}_{33}$,
$u_{14}=-\overline{u}_{31}u_{34}/\overline{u}_{11}$,
$u_{32}=-u_{22}\overline{u}_{23}/\overline{u}_{33}$,
$u_{24}=-\overline{u}_{32}u_{34}/\overline{u}_{22}$. Thus
\beql{}
0=\left\langle C_3,C_4\right\rangle=u_{13}\overline{u}_{14}+u_{23}\overline{u}_{24}+u_{33}\overline{u}_{34}
=-u_{13}\frac{u_{31}\overline{u}_{34}}{u_{11}}-u_{23}\frac{u_{32}\overline{u}_{34}}{u_{22}}+u_{33}\overline{u}_{34}=
\eeq
\begin{equation*}
=\frac{\overline{u}_{34}}{\overline{u}_{33}}\left(\left|u_{13}\right|^2+\left|u_{23}\right|^2+\left|u_{33}\right|^2\right)\neq0,
\end{equation*}

since the last expression in the brackets is strictly positive. This is a contradiction.
\end{proof}

Now we present the dual of Proposition \ref{cond2} and Lemma
\ref{count1}.

\begin{proposition}\label{prop2}
Under the conditions of Theorem \ref{newcond}, a SQ matrix of degree
$n$ satisfies the following properties:

\begin{enumerate}
\item The first row of $X$ is $\left[1,1\right]$;
\item The first row of $Y$ is $\left[0,\hdots, 0\right]$;
\item $n\geq 5$;
\end{enumerate}
\end{proposition}

\begin{proof}
The first two properties are evident from the proof of Proposition
\ref{cond2}. The third one follows from the fact that the only
candidates of order $4$ with these properties are not SQ.
\end{proof}

\begin{lemma}\label{count2}
Suppose that a SQ matrix $M$ of degree $n$ satisfies the conditions
of Theorem \ref{newcond}. Then, its possible number of ones is at
most $n^2-2n+4$, and hence, it has at least $2n-4$ zeros.
\end{lemma}

\begin{proof}
We count the number of ones in each block of $M$ separately. First,
the number of ones in $Q$ is $4$. Then the number of ones in
$J_{3\times 2}$ is $6$. The first condition of Theorem \ref{newcond}
and Proposition \ref{prop2} imply that the number of ones in $X$ and
$Y$ cannot be more than $2$ in each row. Since there are $n-3$ rows
in $X$ and $Y$, we conclude that the possible number of ones is at
most $4+6+2(n-3)+n(n-4)=n^2-2n+4$.
\end{proof}

\begin{corollary}
The following two symmetric, SQ matrices of degree $5$ do not
support unitaries:

{\normalsize
\beql{}
\left\{ \left[
\begin{array}{cc|cc|c}
 1 & 0 & 1 & 1 & 1 \\
 0 & 1 & 1 & 1 & 1 \\
 1 & 1 & 1 & 1 & 0 \\
 \hline
 1 & 1 & 0 & 0 & 1 \\
 1 & 1 & 0 & 0 & 1
\end{array}
\right], \left[
\begin{array}{cc|cc|c}
 1 & 0 & 1 & 1 & 1 \\
 0 & 1 & 1 & 1 & 1 \\
 1 & 1 & 1 & 1 & 1 \\
 \hline
 1 & 1 & 0 & 0 & 1 \\
 1 & 1 & 0 & 0 & 1
\end{array}
\right] \right\}. \eeq }
\end{corollary}

These matrices are equivalent to the exceptional matrices of order $5$ appearing in Section
\ref{list}.

We conclude this section with a SQ matrix of degree $10$ which
satisfies the conditions in both Theorem \ref{cond} and Theorem
\ref{newcond}:

\beql{}
\left[
\begin{array}{cccccccccc}
 1 & 0 & 1 & 1 & 1 & 1 & 0 & 0 & 0 & 0 \\
 0 & 1 & 1 & 1 & 0 & 0 & 1 & 1 & 0 & 0 \\
 1 & 1 & 1 & 1 & 0 & 0 & 0 & 0 & 1 & 1 \\
 0 & 0 & 1 & 0 & 1 & 1 & 1 & 1 & 1 & 1 \\
 0 & 0 & 0 & 1 & 1 & 1 & 1 & 1 & 1 & 1 \\
 1 & 1 & 0 & 0 & 1 & 1 & 1 & 1 & 1 & 1 \\
 1 & 1 & 0 & 0 & 1 & 1 & 1 & 1 & 1 & 1 \\
 1 & 1 & 0 & 0 & 1 & 1 & 1 & 1 & 1 & 1 \\
 1 & 1 & 0 & 0 & 1 & 1 & 1 & 1 & 1 & 1 \\
 1 & 1 & 0 & 0 & 1 & 1 & 1 & 1 & 1 & 1 \\
\end{array}
\right]. \eeq

\end{document}